\documentclass[twocolumn,twoside]{IEEEtran}
\usepackage{amsmath,amssymb,amsthm,epsfig,color,subfigure,empheq,graphicx,comment}
\usepackage{enumerate,url,algorithm,wasysym,epstopdf,balance,enumerate}

\DeclareMathOperator{\rank}{rank}

\DeclareMathOperator{\trace}{Tr}
\DeclareMathOperator{\vectorize}{vec}

\newtheorem{remark}{Remark}

\newcommand \bzero{\mathbf{0}}
\newcommand \bone{\mathbf{1}}
\newcommand \ba{\mathbf{a}}
\newcommand \bb{\mathbf{b}}

\newcommand \bg{\mathbf{g}}

\newcommand \bi{\mathbf{i}}

\newcommand \bdell{\boldsymbol{\ell}} 

\newcommand \bn{\mathbf{n}}

\newcommand \bp{\mathbf{p}}
\newcommand \bq{\mathbf{q}}
\newcommand \br{\mathbf{r}}

\newcommand \bu{\mathbf{u}}
\newcommand \bv{\mathbf{v}}
\newcommand \bw{\mathbf{w}}
\newcommand \bx{\mathbf{x}}
\newcommand \by{\mathbf{y}}

\newcommand \bA{\mathbf{A}}
\newcommand \bB{\mathbf{B}}
\newcommand \bC{\mathbf{C}}

\newcommand \bG{\mathbf{G}}
\newcommand \bH{\mathbf{H}}
\newcommand \bI{\mathbf{I}}
\newcommand \bJ{\mathbf{J}}

\newcommand \bL{\mathbf{L}}

\newcommand \bN{\mathbf{N}}

\newcommand \bP{\mathbf{P}}
\newcommand \bQ{\mathbf{Q}}
\newcommand \bR{\mathbf{R}}
\newcommand \bS{\mathbf{S}}

\newcommand \bU{\mathbf{U}}
\newcommand \bV{\mathbf{V}}

\newcommand \bX{\mathbf{X}}
\newcommand \bY{\mathbf{Y}}
\newcommand \bZ{\mathbf{Z}}

\newcommand \btheta{\boldsymbol{\theta}}

\newcommand \brho{\boldsymbol{\rho}}

\newcommand \bphi{\boldsymbol{\phi}}

\newcommand \bSigma{\mathbf{\Sigma}}

\newcommand \mcD{\mathcal{D}}

\newcommand \mcG{\mathcal{G}}

\newcommand \mcL{\mathcal{L}}
\newcommand \mcM{\mathcal{M}}
\newcommand \mcN{\mathcal{N}}

\newcommand \mcP{\mathcal{P}}

\newcommand \mcS{\mathcal{S}}

\newcommand \bmcS{\bar{\mathcal{S}}}

\newcommand \tbp{\tilde{\mathbf{p}}}

\newcommand \tbN{\tilde{\mathbf{N}}}

\newcommand \tbP{\tilde{\mathbf{P}}}
\newcommand \tbQ{\tilde{\mathbf{Q}}}

\newcommand \tbV{\tilde{\mathbf{V}}}

\newcommand \hbG{\hat{\mathbf{G}}}

\newcommand \hbY{\hat{\mathbf{Y}}}

\newcommand \hbSigma{\hat{\mathbf{\Sigma}}}

\newcommand \bbn{\bar{\mathbf{n}}}

\newcommand \bbp{\bar{\mathbf{p}}}

\newcommand \bbA{\bar{\mathbf{A}}}

\newcommand \bbY{\bar{\mathbf{Y}}}

\begin{document}

\title{Learning Distribution Grid Topologies: A Tutorial}

\author{
	Deepjyoti Deka,~\IEEEmembership{Senior Member,~IEEE}, 
	Vassilis Kekatos,~\IEEEmembership{Senior Member,~IEEE},\\
 and
 Guido Cavraro,~\IEEEmembership{Member,~IEEE}

\thanks{D. Deka is with the Theoretical Division at Los Alamos National Laboratory (LANL), Los Alamos, NM, USA (Email: {\tt\small deepjyoti@lanl.gov}). V. Kekatos is with the Bradley Dept. of ECE, Virginia Tech, Blacksburg, VA, USA (Email: {\tt\small kekatos@vt.edu}). G. Cavraro is with the National Renewable Energy Laboratory (NREL), Golden, CO, USA (Email: {\tt\small guido.cavraro@nrel.gov}). All authors contributed equally to this work.}
}	
	

\maketitle

\begin{abstract}
Unveiling feeder topologies from data is of paramount importance to advance situational awareness and proper utilization of smart resources in power distribution grids. This tutorial summarizes, contrasts, and establishes useful links between recent works on topology identification and detection schemes that have been proposed for power distribution grids. The primary focus is to highlight methods that overcome the limited availability of measurement devices in distribution grids, while enhancing topology estimates using conservation laws of power-flow physics and structural properties of feeders. Grid data from phasor measurement units or smart meters can be collected either passively in the traditional way, or actively, upon actuating grid resources and measuring the feeder's voltage response. Analytical claims on feeder identifiability and detectability are reviewed under disparate meter placement scenarios. Such topology learning claims can be attained exactly or approximately so via algorithmic solutions with various levels of computational complexity, ranging from least-squares fits to convex optimization problems, and from polynomial-time searches over graphs to mixed-integer programs. Although the emphasis is on radial single-phase feeders, extensions to meshed and/or multiphase circuits are sometimes possible and discussed. This tutorial aspires to provide researchers and engineers with knowledge of the current state-of-the-art in tractable distribution grid learning and insights into future directions of work. 
\end{abstract}


\section{Introduction}\label{sec:intro}
Distribution grids constitute the final tier in the delivery of electricity to end-users. To ease protection and voltage control, most distribution grids are operated in a radial (tree-like) topology, which can be modified by changing breaker statuses on available lines~\cite{Kersting}. In recent years, the growth of behind-the-meter distributed energy resources (DERs) and smart loads (e.g., air-conditioners, storage devices, electric vehicles) have brought distribution grids to the forefront of smart grid advancement \cite{stewart2017integrated}. Industrial and academic research on smart distribution grids has advocated the participation of distribution grid resources in wholesale electricity markets and ancillary services (such as demand response, frequency regulation, and transactive energy services~\cite{stewart_liu}). Integrating renewables introduces new challenges for voltage regulation and calls for dispatching DERs without violating physical and operational grid ratings. This necessitates knowing the correct feeder models. Moreover, situational awareness requires accurate distribution grid state estimation (DSSE)~\cite{wang2018survey}, in which the operational topology is a critical component. Topology estimates are also important for ensuring the dynamic stability of inverter-interfaced DERs.

However, distribution utilities often have only partial knowledge of their primary and/or secondary networks and the associated line impedances. Similarly, even if the utility knows the line infrastructure and line impedances, it may not have information on which lines are currently energized. This is owing to the fact that distribution grids are frequently reconfigured for maintenance; load balancing; to improve voltage profiles, minimize losses, or alleviate faults; or rashly, while restoring service after extreme weather events. Such changes may not be logged into the distribution management system, and hence, need to be estimated or at least verified. In this context, grid topology learning can be broadly classified into \emph{topology detection} and \emph{topology identification}. In topology detection, the estimator or system operator knows the line infrastructure and their impedances and needs to determine the ones that are currently energized. Topology identification, on the other hand, aims to determine both the connections and line impedances. Without knowledge of line infrastructure, identification is a more challenging task. 

Reliable topology estimation and DSSE in general are hindered by the limited placement of real-time metering devices on low-voltage grids. This is in stark contrast to high-voltage transmission systems that usually enjoy full observability and breaker statuses on lines are reported in real-time or identified jointly while estimating system states via the so-termed generalized power system state estimators; see \cite{Kekatos_NAPS_2012}, \cite[Ch.~4]{ExpConCanBook}, 
On the other hand, smart meters and real-time sensors as part of the \emph{advanced metering infrastructure} (AMI) have seen increased adoption at customer/load locations. In 2020, about 88\% of the 102.9 million AMI installations in the United States have been at residential customers~\cite{EIA}. Such residential data at feeder buses, though non-ubiquitous, can be used for topology learning and other distribution grid tasks \cite{mohsenian2021smart,arghandeh2017big}.

This tutorial systematically reviews approaches for topology learning in power distribution grids. Emphasizing on topology identification and topology detection, approaches are organized based on the type of measurements needed (namely smart meter data and synchrophasor data), and depending on whether injection and/or voltage readings are collected at all or a subset of buses. The focus is on topology estimates that are provably correct given a sufficient number of measurements. In the realistic setting of partial observability, several of the surveyed algorithms use constraints from the physics of power flows as well as topological properties emerging from the tree structure of feeders to ensure consistency of topology estimates. Depending on meter deployments, oftentimes complete topology recovery is not possible. Nonetheless, several of the presented schemes are able to estimate a reduced feeder model conveying some of the properties of the original grid, which may be sufficient in DER control applications. Identifiability and detectability claims are reviewed to explain when a topology can be unveiled given a particular meter placement. Such claims are important for understanding the limits in topology processing as well as informing future meter placement decisions. 

Topology learning approaches are presented in the context of radial single-phase grids, although generalizations to multiphase and meshed settings are sometimes possible and discussed briefly in Section~\ref{sec:extensions}. The practical relevance of working with single-phase feeder models can be supported by two arguments: \emph{i)} In countries where customers are connected to a three-phase main supply, distribution grids are relatively balanced, and working with the single-phase equivalent model may be a reasonable approximation; \emph{ii)} In countries where residential customers are connected to a single-phase supply, the phase for each customer can be identified relatively easily by clustering voltage data~\cite{blakely2019spectral,wang2016phase,olivier2018phase, bariya2021guaranteed}. Subsequently, per-phase data can be processed assuming the approximate model of a single-phase feeder.

With advances in machine learning, there has been an influx of \emph{supervised} topology learning schemes that train a classifier to detect one out of a finite number of topology configurations. Such classifiers are typically trained using labeled grid data (e.g., voltages) with labels being the associated topologies. Training data can be collected in the field, or synthetically generated using feeder models and simulated loads and renewables. Such approaches are the focus of a recent survey paper~\cite{PLOSsurvey}. In contrast, our tutorial article focuses primarily on physics-aware and unsupervised (i.e., without labeled data) algorithmic solutions for topology learning. 

\emph{Outline:} The rest of this tutorial is organized as follows. Section~\ref{sec:model} introduces some graph and feeder modeling preliminaries. Topology identification schemes relying on synchrophasor and non-synchrophasor data are presented in Sections~\ref{sec:idsyncdata} and \ref{sec:idsmdata}, respectively. Topology detection schemes are reviewed in Section~\ref{sec:detection}. Extensions including generalizations to non-radial topologies, multi-phase grids, data from multiple feeders, and learning from dynamic grid data or non-electric grid data are discussed in Section~\ref{sec:extensions}. Section \ref{sec:topics} discusses research areas related to topology learning and possible directions of future work. The tutorial is concluded in Section~\ref{sec:conclusion}, 

\emph{Notation:} Lower- (upper-) case boldface letters denote column vectors (matrices). Calligraphic symbols are reserved for sets. Symbols $^{\top}$ and $^H$ stands for transposition and complex conjugate transposition. Vectors $\mathbf{0}$ and $\mathbf{1}$ are the all-zero and all-one vectors, while $\bI_{N}$ denotes the $N \times N$ identity matrix. Symbol $\|\mathbf{x}\|$ denotes the Euclidean norm of the vector $\mathbf{x}$, and $\|\mathbf{A}\|_F$ denotes the Frobenius norm of matrix $\bA$. Symbol $\mcD_\bx$ denotes a diagonal matrix with the entries of vector $\bx$ on its main diagonal. The expectation operator is defined as $\mathbb E[\cdot]$. Operator $\vectorize(\bX)$ stacks the columns of matrix $\bX$ into a vector.

\section{Modeling Preliminaries}\label{sec:model}

\subsection{Preliminaries on Graphs}
Because graph algorithms have been widely used in the literature to reconstruct grid topologies, this section reviews some basic concepts of graph theory; see Figure~\ref{fig:graph1}. Consider graph $\mcG_o=(\mcN,\mcL)$. A rooted tree is a connected graph without loops with one node designated as the root and indexed by 0. In a tree graph, a \emph{path} $\mcP_{n,m}$ is the unique sequence of edges connecting any two nodes $n$ and $m$. If node $n$ is along path $\mcP_{0,m}$, then node $m$ is a \emph{descendant} of $n$, and $n$ is an \emph{ancestor} of $m$. If there exists line $(m,n)\in \mcL$ and node $m$ is closer to the root, then node $m$ is the \emph{parent} of $n$, and $n$ is the \emph{child} of $m$. In radial grids, each node modulo the root has a unique parent. A node without children is called a \emph{leaf}, while non-leaf nodes will be termed \emph{internal} nodes. 

\subsection{Feeder Modeling}\label{subsec:feeder}
A single-phase power distribution grid having $N+1$ buses can be modeled by a graph $\mcG_o=(\mcN,\mcL)$. The nodes in $\mcN:=\{0,\ldots,N\}$ correspond to buses, and the edges in $\mcL$ to distribution lines. Let $L$ be the cardinality of $\mcL$; for radial networks $L=N$. The substation is labeled as node $0$. Let $u_n=v_n e^{j \theta_n}$ be the voltage phasor (polar form), and $i_n$ and $s_n=p_n+jq_n$ be respectively, the complex current and complex power injected from bus $n$ into the grid for all $n\in\mcN$. Moreover, let $z_\ell = r_\ell+j x_\ell$ be the impedance and $y_\ell=g_\ell+j b_\ell$ the admittance of line $\ell\in\mcL$. 

Grid topology is captured by the branch-bus incidence matrix $\bbA \in\{0,\pm1\}^{L\times (N+1)}$, which can be partitioned into its first and the rest of its columns as $\bbA=[\mathbf{a}_0~\bA]$. For a radial grid $(L=N)$, the \emph{reduced incidence matrix} $\bA$ is square and invertible~\cite{godsil_2001_algebraic}. Because $\bbA \bone = \bzero$, it holds that $\ba_0 = -\bA^{-1} \bone$. 

\begin{figure*}[t]
\centering	
\includegraphics[width=0.3\columnwidth]{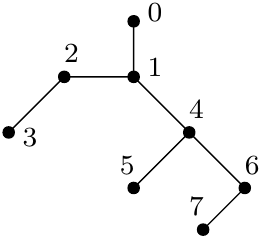}
\hspace*{1em}
\includegraphics[width=0.65\columnwidth]{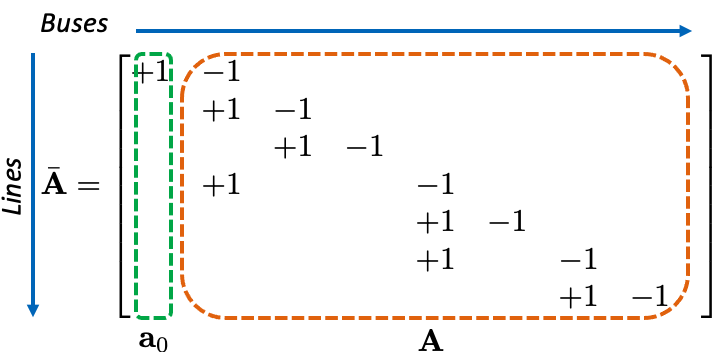}
\hspace*{1em}
\includegraphics[width=0.45\columnwidth]{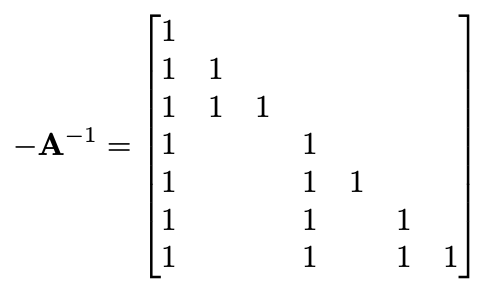}
\hspace*{1em}
\includegraphics[width=0.4\columnwidth]{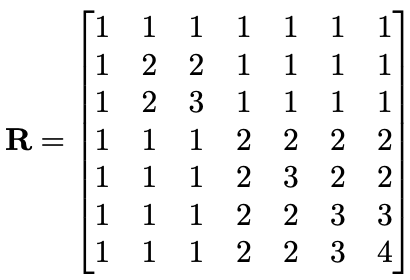}
\caption{\emph{Left to right:} \emph{i)} A graphical depiction of a radial grid with $(3,5,7)$ as leaf nodes and 0 as the root node. \emph{ii)} Its complete and reduced branch-bus incidence matrices. Each row corresponds to an edge and has entry +1 at the parent and -1 at the child node. \emph{iii)} The reduced incidence matrix inverted: The ones on row $n$ mark the ancestors of bus $n$; the ones on column $n$ mark the descendants of bus $n$. \emph{iv)} Matrix $\bR$ assuming all lines have unit resistances: Entry $\bR(5,7)=2$ because paths $\mcP_{0,5}$ (from substation to bus 5) and $\mcP_{0,7}$ (from substation to bus 7) share lines $(0,1)$ and $(1,4)$ of total resistance 2.}
\label{fig:graph1}
\end{figure*}

Voltage and current phasors are related through the bus admittance matrix $\bbY:=\bbA^\top\mcD_\by\bbA$, where $\mcD_\by$ is a diagonal matrix carrying the vector of line admittances $\by=\bg+j\bb$ on its main diagonal; shunt admittances have been ignored. Upon partitioning to separate the substation, we get that
\begin{equation}\label{eq:i=Yv0}
\begin{bmatrix}
i_0\\
\bi
\end{bmatrix}=\underbrace{
\begin{bmatrix}
y_0 & \by_0^\top\\
\by_0 & \bY
\end{bmatrix}}
_{\bbY}
\begin{bmatrix}
u_0\\
\bu
\end{bmatrix}\quad \Rightarrow \quad \bi=\bY(\bu-u_0\bone)
\end{equation}
where $\bY=\bA^\top\mcD_\by\bA$ is the reduced bus admittance matrix, and vectors $(\bi,\bu)$ collect current and voltage phasors at all buses but the substation. The right-hand side (RHS) of \eqref{eq:i=Yv0} follows from $\by_0=-\bY\bone$. We will oftentimes assume $u_0=1$ or that $u_0$ is known and can thus be subtracted from $\bu$. Note that ignoring the substation current in \eqref{eq:i=Yv0} does not incur any loss of information as $i_0=-\bone^\top\bi$. 
Power injections depend nonlinearly on voltages as $s_n=u_ni_n^*$. To arrive at simpler linear models, we usually linearize the nonlinear AC power flow equations. If linearization occurs at the so-termed \emph{flat voltage profile} where $v_n=1$ and $\theta_n=0$ for all $n$, it can be shown that~\cite{Bolognani_2015_Allerton,Deka_TCONES_2018,Bolognani_2015_TAC}
\begin{subequations}\label{eq:lpf}
\begin{align*}
\bp&\simeq \bA^\top\mcD_\bg\bA (\bv-v_0\bone) + \bA^\top\mcD_\bb\bA\btheta\\
\bq&\simeq \bA^\top\mcD_\bb\bA (\bv-v_0\bone) - \bA^\top\mcD_\bg\bA\btheta
\end{align*}
\end{subequations}
where vectors $(\bp,\bq,\bv,\btheta)$ collect power injections, voltage magnitudes and angles at all buses excluding the substation. 

The aforesaid model relates power injections linearly to voltages and holds for radial and meshed single-phase grids alike. The model can be inverted so that voltages are expressed as linear functions of power injections as
\begin{subequations}\label{eq:ldf0}
\begin{align}
\bv&\simeq \bR \bp + \bX \bq + v_0\bone\label{eq:ldf0:v} \\
\btheta&\simeq \bX \bp - \bR \bq\label{eq:ldf0:theta}
\end{align}
\end{subequations}
for appropriately defined matrices $(\bR,\bX)$; see~\cite{Bolognani_2015_Allerton,Deka_TCONES_2018}.\footnote{This linearization is related to the linearization used to derive the so-termed \emph{DC power flow model} in transmission systems. However, the DC power flow model relies on the additional assumption of small resistance-to-reactance ratios, under which matrix $\bR$ can be dropped and \eqref{eq:ldf0:theta} yields the widely used $\bp=\bB\btheta$ model with $\bB=\bX^{-1}$. For distribution systems, both matrices $(\bR,\bX)$ are important.} The substation voltage is usually assumed to be $v_0=1$, or fixed at a known value. 

For the special yet practically relevant case of \emph{radial single-phase grids}, matrices $(\bR,\bX)$ in \eqref{eq:ldf0} can be expressed as
\begin{equation}\label{eq:R+X}
\bR=\bA^{-1}\mcD_\br\bA^{-\top}\quad\text{and}\quad
\bX=\bA^{-1}\mcD_\bx\bA^{-\top}
\end{equation}
where $(\mcD_\br,\mcD_\bx)$ carry line resistances and reactances, respectively, on their main diagonal. The majority of the methods to be discussed in this tutorial refer to radial single-phase feeders for which \eqref{eq:R+X} holds. Due to the properties of $\bA^{-1}$ (see Figure~\ref{fig:graph1}), the $(m,n)$-th entry of $\bR$ enjoys the analytical form~\cite{Deka_TCONES_2018}:
\begin{equation}\label{R_add}
\bR(m,n)= \sum_{(k,l) \in {\mathcal P}_{m,0}\bigcap {\mathcal P}_{n,0}} r_{kl}, \qquad \forall m,n\in\mcN.
\end{equation}
A similar expression holds for $\bX$. This property is exploited later to derive greedy graph algorithms for topology identification. Apparently, matrices $(\bR,\bX)$ have non-negative entries for radial grids. This implies that both active and reactive power injections have a non-negative effect on voltage magnitudes based on \eqref{eq:ldf0}. It also explains in part why voltage angles are small in distribution grids that typically host mostly inductive loads. 

Because of their importance in topology learning, let us also define explicitly the inverses of matrices $(\bR,\bX)$:
\begin{subequations}\label{eq:G+B}
\begin{align}
\bG&:=\bR^{-1}=\bA^\top\mcD_\br^{-1}\bA\label{eq:G+B:G}\\
\bB&:=\bX^{-1}=\bA^\top\mcD_\bx^{-1}\bA.\label{eq:G+B:B}
\end{align}
\end{subequations}
Matrices $(\bG,\bB)$ are positive definite. They can be thought of as reduced Laplacian matrices of the grid graph $\mcG_o$ with edges weighted by the inverse line resistances or reactances. 

Before concluding this section, it is worth contrasting the linearized model in \eqref{eq:ldf0} with the widely used \emph{linear distribution flow (LDF) or LinDistFlow model}, which is derived from the \emph{distribution flow (DistFlow) model} upon ignoring ohmic losses on lines. Per the DistFlow model~\cite{Farivar_2013_CDC}
\begin{subequations}\label{eq:DF}
\begin{align}
\bp&=\bA^\top\bP +\mcD_\br \bdell\label{eq:DF:p}\\
\bq&=\bA^\top\bQ +\mcD_\bx \bdell\label{eq:DF:q}\\
\bA\bv^2+v_0^2\ba_0 &= 2\mcD_\br\bP + 2\mcD_\bx \bQ -(\mcD_\br^2+\mcD_\bx^2)\bdell\label{eq:DF:v}
\end{align}
\end{subequations}
where notation $\bv^2$ denotes a vector carrying the \emph{squared voltage magnitudes} at all buses excluding the substation; the $n$-th entries of $(\bP,\bQ,\bdell)$ carry respectively the active flow, the reactive flow, and the squared current magnitude of the \emph{line} feeding bus $n$. By dropping losses and substituting $(\bP,\bQ)$ from \eqref{eq:DF:p}--\eqref{eq:DF:q} into \eqref{eq:DF:v}, one obtains the LDF model
\begin{equation}\label{eq:ldf02}
\bv^2\simeq 2\bR \bp + 2\bX \bq + v_0^2\bone
\end{equation}
with $(\bR,\bX)$ defined as in \eqref{eq:R+X}. Contrasted with \eqref{eq:ldf0:v}, LDF involves factors of $2$ as it applies to squared voltage magnitudes. The LDF model can alternatively be derived directly from \eqref{eq:ldf0:v} upon linearizing squared voltages around one pu to get $v_n^2\simeq 2v_n-1$; see~\cite{Bolognani_2015_TAC}. The LDF approximation of squared voltages can be shown to be an overestimator of the actual squared voltages~\cite{Gan12}. Hence, the linearization in \eqref{eq:ldf0:v} is also an overestimator of the actual voltages.

So far, we have explicated model \eqref{eq:ldf0} for radial single-phase feeders. Nonetheless, approximate linearized grid models can also be derived for \emph{meshed single-phase}~\cite{Deka_TCONES_2018},~\cite[Section~IV.A]{Cavraro_TSG_2019}; as well as \emph{radial multiphase feeders}~\cite{VKZG16}.

Having established commonly used feeder models, Sections~\ref{sec:idsyncdata} and \ref{sec:idsmdata} deal with topology identification, Section~\ref{sec:detection} deals with topology detection, while extensions are discussed in Section~\ref{sec:extensions}. Table~\ref{table:1} summarizes the key features of the methods to be presented in this tutorial article.

\begin{table*}[t]
\caption{Key features of the presented feeder topology identification and verification methods.}
\centering
\begin{tabular}{|r|| c| c| c| c| c| c|| c| c| c| c| c| c|| c|} 
\hline \hline
& \multicolumn{6}{c||}{\emph{Data}} & \multicolumn{6}{c||}{\emph{Feeder}} &\\
\hline
& \multicolumn{2}{c}{{Type}} &
\multicolumn{2}{|c|}{{Availability}} &
\multicolumn{2}{c||}{{Collection}} &
\multicolumn{2}{c}{{Topology}} &
\multicolumn{2}{|c|}{{Phases}} &
\multicolumn{2}{|c||}{{Feeder}} & {Recovery}\\
\hline
\emph{Method} & {AMI}& {PMU } & {Full} & {Partial} & {Passive} & {Active} & {Tree} & {Mesh}& {$1\phi$} & {$3\phi$} & {One} & {Many} & claims\\ [0.5ex] 
 \hline 
 \multicolumn{14}{|c|}{\textbf{\emph{Feeder Topology Identification}}} \\ [0.5ex]
 \hline
\cite{bariya2021guaranteed,Deka_TPS_2020,gandluru2019joint,liao2019unbalanced} &\checkmark &\checkmark &\checkmark & &\checkmark & &\checkmark & & &\checkmark &\checkmark & &\checkmark\\
\cite{Srinivas22} &\checkmark &\checkmark & &\checkmark &\checkmark & &\checkmark & &\checkmark & &\checkmark & &\\
\cite{arya,Pengwah2022SJ,7981061-R2} &\checkmark & & &\checkmark &\checkmark & &\checkmark & &\checkmark & &\checkmark & &\\
\cite{Bolognani_CDC_2013} &\checkmark & &\checkmark & &\checkmark & &\checkmark & &\checkmark & &\checkmark & &\checkmark\\
\cite{Cavraro_LCSS_2018} &\checkmark & &\checkmark &\checkmark & &\checkmark &\checkmark & &\checkmark & &\checkmark & &\checkmark\\
\cite{Cavraro_TCONES_2019,Taheri_MILP_2019} &\checkmark & &\checkmark & & &\checkmark &\checkmark & &\checkmark & &\checkmark & &\checkmark\\
\cite{Deka_SmartGridComm_2016,Deka_ECC_2016,Deka_TCONES_2018,Deka_TCONES_2020} --&\checkmark & & &\checkmark &\checkmark & &\checkmark & &\checkmark & & &\checkmark &\checkmark\\
-- \cite{park2018exact,Park_TCONES_2020,park2018learning} & & & & & & & & & & & & &\\
\cite{Deka_TCONES_2021,Deka_TSG_2020,liao2018urban,Liao_PESGM_2016} &\checkmark & &\checkmark & &\checkmark & & &\checkmark &\checkmark & & &\checkmark &\checkmark\\
\cite{Li_JEPES_2021} &\checkmark & & &\checkmark &\checkmark & &\checkmark & & &\checkmark &\checkmark & &\checkmark\\
\cite{Tirja2} &\checkmark & & &\checkmark &\checkmark & & &\checkmark &\checkmark & &\checkmark & &\checkmark\\
\cite{Korres2012GSE,9748966,Zhang_TPS_2021,Zhang_TSG_2020, Zhao_TSG_2020} &\checkmark & &\checkmark & &\checkmark & &\checkmark & &\checkmark & &\checkmark & &\\ 
\cite{Ardakanian_TCONES_2019,Brouillon_CDC_2021,Yu_TPS_2018,Yu_TPS_2019} & &\checkmark &\checkmark & &\checkmark & & &\checkmark & &\checkmark &\checkmark & &\\
\cite{bariya2018data} & &\checkmark &\checkmark & &\checkmark & &\checkmark & &\checkmark & &\checkmark & &\\
\cite{Deka_PSCC_2016,Liao_MAPS_2015,Talukdar_ACC_2017,Weng_TPS_2017,Zhao_TSG_2020} & &\checkmark &\checkmark & &\checkmark & &\checkmark & &\checkmark & & &\checkmark &\checkmark\\ 
 \cite{Fabbiani_TCST_2021,Tirja1,Lateef_ISGT_2019} & &\checkmark &\checkmark & &\checkmark & & &\checkmark &\checkmark & &\checkmark & &\\
 \cite{anguluri2021grid,doddi2020learning,9034147} & &\checkmark & &\checkmark &\checkmark & & &\checkmark &\checkmark & &\checkmark & &\checkmark\\
\cite{Miao_ACC_2019} & &\checkmark & &\checkmark &\checkmark & &\checkmark & &\checkmark & &\checkmark & &\\
\cite{Moffat_TSG_2020,yuan_TCONES_2016} & &\checkmark & &\checkmark &\checkmark & &\checkmark & &\checkmark & &\checkmark & &\checkmark\\
\cite{Talukdar_ACC_2018} & &\checkmark & &\checkmark &\checkmark & &\checkmark & &\checkmark & & &\checkmark &\checkmark\\ 
\cite{Talukdar_Automatica_2020} & &\checkmark &\checkmark & &\checkmark & & &\checkmark &\checkmark & & &\checkmark &\checkmark\\
\hline
 \multicolumn{14}{|c|}{\textbf{\emph{Feeder Topology Verification}}} \\ [0.5ex]
\hline 
\cite{Arghandeh_PESGM_2015,Cavraro_ISGT_2015} & &\checkmark & &\checkmark &\checkmark & & &\checkmark &\checkmark & &\checkmark & &\\
\cite{Cavraro_PESGM_2015,Cavraro_TPS_2018} & &\checkmark & &\checkmark &\checkmark & & &\checkmark &\checkmark & &\checkmark & &\checkmark\\
\cite{Cavraro_LCSS_2020,He2021hybrid,Tian_TPS_2016} &\checkmark & & &\checkmark &\checkmark & &\checkmark & &\checkmark & &\checkmark & &\checkmark\\
\cite{Cavraro_TCONES_2019} &\checkmark & &\checkmark & & &\checkmark &\checkmark & &\checkmark & &\checkmark & &\checkmark\\
\cite{Cavraro_TSG_2019} &\checkmark & &\checkmark & &\checkmark & &\checkmark &\checkmark &\checkmark & &\checkmark &\checkmark &\checkmark\\
\cite{Farajollahi_TSG_2020,9852281} & &\checkmark & &\checkmark &\checkmark & &\checkmark & &\checkmark & &\checkmark & &\\
\cite{Korres2012GSE,9748966,Sharon_ISGT_2012} &\checkmark & & &\checkmark &\checkmark & & &\checkmark &\checkmark & &\checkmark & &\checkmark\\
\cite{sevlian2015distribution} &\checkmark & & &\checkmark &\checkmark & &\checkmark & &\checkmark & & &\checkmark&\checkmark\\
\cite{Soltani2022GSE} &\checkmark & & &\checkmark &\checkmark & & &\checkmark & &\checkmark &\checkmark & &\\
\hline\hline
\end{tabular}
\label{table:1}
\end{table*}

\section{Topology Identification using Phasor Data}\label{sec:idsyncdata}
Given voltage with/without injection data collected from a feeder, topology identification aims at finding the network connectivity with/without line impedances of the feeder. This section reviews methods relying on synchrophasor data, while Section~\ref{sec:idsmdata} reviews methods relying on smart meter data. For synchrophasor data, the data model of interest is $\bi=\bY\bu$ as the operator collects synchrophasor measurements of nodal voltages and/or currents on a subset of measured buses. Identifying the topology amounts to estimating the complex-valued grid Laplacian matrix $\bY$ or at least its sparsity pattern (non-zero entries). Depending on the type and location of synchrophasor data, the ensuing setups may be identified.

\subsection{Setup S1) Complete Data}\label{subsec:S1} 
This ideal setup presumes voltage and current phasors $\{(\bi_t,\bu_t)\}_{t=1}^T$ are collected at all buses. Data collected over $T$ times are stacked as columns of $N\times T$ matrices $\bI$ and $\bU$. To account for measurement noise, we can postulate the model 
\begin{equation}\label{eq:i=Yv}
\bI=\bY\bU+\bN.
\end{equation}
If $\rank(\bU)=N$, matrix $\bY$ can be found as the least-squares estimate $\hbY_{\text{LS}}:=\bI\bU^H\left(\bU\bU^H\right)^{-1}$, which is the minimizer of $\|\bI-\bY\bU\|_F^2$; see~\cite{Lateef_ISGT_2019,9034147}. Guaranteeing $\rank(\bU)=N$ requires $T\geq N$ and that the feeder undergoes sufficiently different loading conditions so that voltage sequences across buses are linearly independent. Instead of batch processing, reference \cite{Fabbiani_TCST_2021} puts forth an adaptive scheme to deal with streaming data and updates $\bY$ in an online fashion using a recursive least squares algorithm. 

Either batch or adaptively computed, the least-squares (LS) estimate $\hbY_{\text{LS}}$ will not have any zero entries due to numerical/noise errors, whereas the actual $\bY$ is known to be sparse. To promote sparse estimates, the LS fit can be regularized by the $\ell_1$-norm of the entries of $\bY$. An $\ell_1$-norm regularizer can be used to find meaningful estimates~\cite{Ardakanian_TCONES_2019}. This is particularly useful when the sample size $T$ is smaller than the number of unknown variables in $\bY$. Additional regularization terms and constraints could be appended to the data fitting problem to enforce or promote the properties of the actual $\bY$. Reference~\cite{Tirja1} aims at estimating $\bY$ using voltage phasors and power rather than current injections based on an LS fit properly regularized and constrained to account for the Laplacian nature of $\bY$. An alternative two-step heuristic to yield sparse estimates of $\bY$ is suggested in~\cite{Yu_TPS_2018}: After estimating $\bY$ using an LS fit, lines whose estimated conductances are relatively small are progressively removed. Heed that \eqref{eq:i=Yv} models measurement noise on current phasors alone, even though both voltage and current readings are noisy. To deal with both sources of noise, a total least-squares (TLS) or \emph{error-in-variables} approach is pursued in~\cite{Yu_TPS_2018,Yu_TPS_2019}. A Bayesian version of TLS with sparsity-promoting priors is suggested in~\cite{Brouillon_CDC_2021,SparseTLS}.

References~\cite{Angjelichinoski_IDCM_2017} and \cite{Angjelichinoski_2018} consider the same topology identification setup for the case of a DC microgrid. To ensure $\bU$ is full row-rank, the grid is intentionally driven to sufficiently different states. The latter is achieved by changing the droop parameters of inverter-interfaced generators. 

\subsection{Setup S2) Partial Injection Data}\label{subsec:S2} 
This setup assumes that voltage data are measured at all buses, but current data are only measured on a subset of buses $\mcS\subset\mcN$ with cardinality $S=|\mcS|$. Let $\bmcS=\mcN\setminus \mcS$ denote its complement set. Then, the data model can be partitioned as
\begin{equation}\label{eq:partition}
\begin{bmatrix}
\bI_\mcS\\
\bI_{\bmcS}
\end{bmatrix}=
\begin{bmatrix}
\bY_{\mcS\mcS} & \bY_{\mcS\bmcS}\\
\bY_{\bmcS\mcS} & \bY_{\bmcS\bmcS}
\end{bmatrix}
\begin{bmatrix}
\bU_\mcS\\
\bU_{\bmcS}
\end{bmatrix}.
\end{equation}
where metered and non-metered buses have been permuted accordingly without losing generality. Given $(\bI_\mcS,\bU_\mcS,\bU_{\bmcS})$, matrices $\bY_{\mcS\mcS}$ and $\bY_{\mcS\bmcS}$ can be estimated using an LS fit.

\subsection{Setup S3) Partial Data}\label{subsec:S3}
This may be the most practically relevant setup. Voltage and current data are collected only at a subset of buses $\mcS$. Most methods further assume that all buses in $\bmcS$ are \emph{zero-injection buses (ZIB)}. Plugging $\bI_{\bmcS}=\bzero$ into \eqref{eq:partition} yields $\bU_{\bmcS}=-\bY_{\bmcS\bmcS}^{-1}\bY_{\bmcS\mcS}\bU_\mcS$.\footnote{Because of this, if non-ZIB voltages $\bU_{\mcS}$ are metered, measuring ZIB voltages $\bU_{\bmcS}$ does not really offer any additional information.} Substituting $\bU_{\bmcS}$ into \eqref{eq:partition} provides the Kron-reduced data model~\cite{Dorfler_2013_TCS}
\begin{equation}\label{eq:Kron}
\bI_\mcS=\bbY_{\mcS}\bU_\mcS
\end{equation}
where $\bbY_{\mcS}:=\bY_{\mcS\mcS}-\bY_{\mcS\bmcS}\bY_{\bmcS\bmcS}^{-1}\bY_{\bmcS\mcS}$. The \emph{Kron-reduced or effective} admittance matrix $\bbY_{\mcS}$ is still a (reduced) weighted Laplacian of a graph~\cite{Moffat_TSG_2020}. Nonetheless, this graph is in general non-radial despite the complete graph being radial. Topology recovery algorithms typically proceed in two steps. Matrix $\bbY_\mcS$ is first estimated using an LS fit on data $(\bI_\mcS,\bU_\mcS)$. Given this estimate, graph algorithms are then used to infer a radial topology~\cite{Moffat_TSG_2020,yuan_TCONES_2016}. Reference \cite{Moffat_TSG_2020} proposes a complex-valued extension of the \emph{recursive grouping} algorithm, detailed later in Section~\ref{subsec:m3}. The approach suggested in~\cite{Miao_ACC_2019} aims at recovering the entire topology under the additional assumption that the number of ZIBs is known.

\begin{remark}\label{re:tree}
As noted earlier, a general matrix $\bY$ can be recovered from \eqref{eq:i=Yv} if and only if $\rank(\bU)=N$. Unfortunately, such an assumption cannot be satisfied in the presence of ZIBs, as ZIB voltages are linearly dependent on non-ZIB voltages. If there are say $K$ ZIBs, it holds that $\rank(\bU)\leq N-K$, and hence, the LS estimate $\hbY_\text{LS}$ is undefined even if the complete voltage and current data are available. However, not all is lost. The previous identifiability analysis treated $\bY$ as a general matrix. If one exploits the fact that $\bY$ is the Laplacian of a tree graph, there exist identifiability claims and algorithms that can recover $\bY$ using partial data. For example, reference~\cite{yuan_TCONES_2016} can unveil the actual radial topology from $\bbY_\mcS$ if each ZIB has at least three neighbors. References~\cite{Cavraro_TCONES_2019,park2018exact,Park_TCONES_2020} provide a weaker assumption as they can recover the complete $\bY$ if all leaf nodes are non-ZIBs. The aforesaid references rely on a smart meter data model and they are presented in more detail in Section~\ref{subsec:m3}, yet they can be properly extended to the synchrophasor data model along the lines of \cite{Moffat_TSG_2020}.
\end{remark}

\section{Topology Identification via Smart Meter Data}\label{sec:idsmdata}
Even though synchrophasor data have desirable monitoring features (high sampling frequency, accuracy, and time synchronization), they come at increased cost so widespread deployments are currently rare. Smart meters on the other hand are more widely available. A smart meter may be measuring active power injections, reactive power injections, and/or voltage magnitudes. The data model of choice is the linearized grid model \eqref{eq:ldf0:v}, which we repeat for time $t$ as
\begin{equation}\label{eq:ldf}
\bv_t=\bR\bp_t+\bX\bq_t+\bone+\bn_t,
\end{equation}
assuming the substation voltage is fixed at $v_0=1$~pu without loss of generality. The noise term $\bn_t$ captures measurement noise and modeling (e.g., linearization and single-phase approximation) errors. Sample data across times $t=1,\ldots,T$, and stack them as columns of the $N\times T$ matrices $(\bV,\bP,\bQ)$ to arrive at
\begin{equation}\label{eq:ldf2}
\bV=\bR\bP+\bX\bQ+\bone\bone^\top+\bN.
\end{equation}
Although the emphasis is on smart meter data, model \eqref{eq:ldf2} may be relevant to any non-synchronized grid data. Reference~\cite{Srinivas22} arrives at a model similar to \eqref{eq:ldf2} upon successive linearization of the power flow equations and estimates the topology matrices using an unscented Kalman filter-based approach.

Different from the phasor data model, model~\eqref{eq:ldf2} involves real-valued quantities. A key challenge is that one has to estimate two topology matrices $(\bR,\bX)$ from $NT$ equations. To overcome this issue, two types of approaches (grid probing and covariance-based methods) have been developed and will be reviewed later on. Similar to topology identification using phasor data, the ensuing setups can occur.

\subsection{Setup M1) Complete Data}\label{subsec:m1}
Similar to setup~\emph{S1)} in Section \ref{subsec:S1}, suppose one collects power and voltage data at all buses so that $(\bV,\bP,\bQ)$ are known. Suppose we focus on finding $\bR$. To eliminate the effect of reactive power on voltage data, one could try the technique of \emph{grid probing}~\cite{Cavraro_TCONES_2019,Cavraro_LCSS_2018}. Grid probing works with a data model on \emph{successive time differences}
\begin{equation}\label{eq:dldf}
\tbV=\bR\tbP+\bX\tbQ+\tbN
\end{equation}
where the $t$-th column of $\tbP$ is $\tbp_t:=\bp_{t}-\bp_t$. In other words, matrix $\tbP$ contains the differences in active power injections between consecutive times. Matrices $\tbV$ and $\tbQ$ are defined similarly. Note that time differentiation amplifies (doubles) the variance in the entries of $\tbN$ compared to $\bN$.
Suppose for now that the operator can control the injections at all buses. If reactive powers are kept or assumed constant across time, model~\eqref{eq:dldf} simplifies as 
\begin{equation}\label{eq:dldf2}
\tbV=\bR\tbP+\tbN~~~\Leftrightarrow~~~\tbP=\bG\tbV-\bG\tbN
\end{equation}
which upon inversion yields the real-valued counterpart of \eqref{eq:i=Yv}. In theory, any method discussed under \emph{S1)} is applicable here. In practice, the noise term $\tbN$ can be stronger than the noise term in \eqref{eq:i=Yv}, and $\bG\tbN$ is now correlated across buses. 

An alternative idea to eliminate the effect of reactive power injections is a covariance-based approach. All grid data (voltages and powers) have been \emph{centered}, i.e., for each signal, its mean (time-average) has been subtracted. Due to this preprocessing, the correlation between centered data coincides with the covariance between the original data. All processing is henceforth performed on centered data. Let us compute the correlation between voltage and power injection vectors. From~\eqref{eq:ldf}, these matrices can be readily computed as
\begin{subequations}\label{eq:corr}
\begin{align}
\bSigma_{vp}&=\bR\bSigma_{pp}+\bX\bSigma_{pq}^\top\label{eq:corr:a}\\
\bSigma_{vq}&=\bR\bSigma_{pq}+\bX\bSigma_{qq}.\label{eq:corr:b}
\end{align}
\end{subequations}
where $\bSigma_{vp}:=\mathbb{E}[\bv_t\bp_t^\top]$ can be practically estimated from data as $\hbSigma_{vp}=\frac{1}{T}\sum_{t=1}^T\bv_t\bp_t^\top$. To simplify notation, we will henceforth denote $\hbSigma_{vp}$ simply by $\bSigma_{vp}$. The remaining correlation matrices $(\bSigma_{vq},\bSigma_{pq},\bSigma_{qq})$ can be estimated from data similarly. In the absence of solar generation, it is reasonable to assume that matrices $\bSigma_{pp}$, $\bSigma_{qq}$, and $\bSigma_{pq}$ are diagonal. Then, it is not hard to verify from \eqref{eq:corr} that the $(m,n)$-th entries of matrices $\bSigma_{vp}$ and $\bSigma_{vq}$ satisfy the linear equations
\begin{equation}\label{eq:corrmn}
\begin{bmatrix}
\mathbb{E}[v_mp_n]\\
\mathbb{E}[v_mq_n]
\end{bmatrix}=
\begin{bmatrix}
\mathbb{E}[p_n^2] & \mathbb{E}[p_nq_n]\\
\mathbb{E}[p_nq_n] & \mathbb{E}[q_n^2]
\end{bmatrix}
\begin{bmatrix}
\bR(m,n)\\
\bX(m,n)
\end{bmatrix}.
\end{equation}
If $\mathbb{E}[p_n^2]\cdot \mathbb{E}[q_n^2]>(\mathbb{E}[p_nq_n])^2$, the matrix in \eqref{eq:corrmn} is invertible, and we can readily solve for $\bR(m,n)$ and $\bX(m,n)$ and then derive the topology by matrix inversion. 

References \cite{Zhang_TSG_2020,Zhang_TPS_2021} rely also on complete data $(\bV,\bP,\bQ)$ to recover grid topology using a two-step approach involving a linearization wherein ratios $p_n/v_n$ and $q_n/v_n$ are approximated as linear functions of voltage magnitudes and angles.

\subsection{Setup M2) Partial Injection Data}\label{subsec:m2}
Paralleling \emph{S2)} in Section \ref{subsec:S2}, this setup assumes voltage data are collected at all buses, whereas power injection data are collected only at a subset of buses $\mcS$. The linearized models in \eqref{eq:ldf}--\eqref{eq:dldf2} can be partitioned similarly to \eqref{eq:partition}. Therefore, for setup \emph{M2)} the given data is $(\tbV,\bP_\mcS,\bQ_\mcS)$. Based on these data, it is not hard to see that the covariance-based approach can estimate matrices 
\begin{equation}\label{eq:RsXs}
\bR_\mcS:=\begin{bmatrix}
\bR_{\mcS\mcS}\\
\bR_{\bmcS\mcS}
\end{bmatrix}~~\text{and}~~
\bX_\mcS:=\begin{bmatrix}
\bX_{\mcS\mcS}\\
\bX_{\bmcS\mcS}
\end{bmatrix}.
\end{equation}

Setup \emph{M2)} may be more natural for the grid probing approach since it is unlikely an operator can control injections at every single bus. In that case, the set $\mcS$ is the set of buses that the operator can actually control. Power injections at $\mcS$ can be intentionally altered for short intervals. Building again on the difference model and assuming the operator perturbs only active injections at the buses in $\mcS$, we obtain
\begin{equation}\label{eq:m2}
\tbV=\bR_\mcS\tbP_\mcS+\tbN
\end{equation}
The columns of $\bR$ associated with $\bmcS$ have been suppressed assuming injections in $\mcS$ are perturbed for short intervals during which injections in $\bmcS$ remain roughly unchanged. 

Despite the partial injection data, the complete matrix $\bG=\bR^{-1}$ is identifiable if $\rank(\tbP_\mcS)=S$ and the grid is probed at all leaf buses~\cite{Cavraro_TCONES_2019}. This identifiability result carries over to the synchrophasor data model. Note that, thanks to the link established between setups \emph{M2)} and \emph{S2)} in Remark~\ref{re:tree}, the algorithms to be presented next are also applicable to \emph{M2)}. 

Albeit the previous result asserts that the complete grid topology can be identified, it does not explain how. We next outline three alternative approaches for topology recovery under setup \emph{M2)}. The first group of approaches consists of the graph algorithms put forth in~\cite{park2018exact,Park_TCONES_2020,Cavraro_LCSS_2018}. These algorithms use the knowledge of entries in $\bR_{\mcS\mcS}$ to learn the complete $\bR$ leveraging properties of radial networks via a polynomial-complexity algorithm. These graph algorithms are discussed later under setup \emph{M3)}.

Despite their polynomial complexity, graph algorithms can be sensitive to noisy entries of $\bR_\mcS$. As an alternative, one can try recovering $\bG=\bR^{-1}$ directly by fitting the available data $(\tbV,\tbP_\mcS)$; see~\cite{Cavraro_TCONES_2019}: Premultiplying \eqref{eq:m2} by $\bG$ yields a model similar to the RHS of \eqref{eq:dldf2} with $\tbP^\top=[\tbP_\mcS^\top~
\bzero^\top]$, which can be used to find an LS fit for $\bG$. Nonetheless, the estimate $\hbG_\text{LS}$ would not satisfy the properties of the actual $\bG$. In fact, ensuring the estimate $\hbG$ is the reduced Laplacian of a tree graph would entail a non-convex feasible set. To arrive at a convex data fitting problem, reference~\cite{Cavraro_TCONES_2019} suggests finding $\bG$ as the minimizer of
\begin{equation}\label{eq:LaplacianEstimation}
\min_{\bG\in\mcM}~\|\tbP-\bG\tbV\|_F^2 + \lambda\trace[\bG(\bI+\bone\bone^\top)] - \mu \log|\bG|
\end{equation}
where $\lambda,~\mu>0$ are tunable parameters. Set $\mcM$ ensures the minimizer is symmetric with non-positive off-diagonal entries. It can also incorporate prior information, such as $\bG_{mn}=0$ if buses $(m,n)$ are known not to be connected. The second term in the objective sums up the absolute off-diagonal entries of $\bG$ to promote sparsity. The third term serves as a barrier keeping $\bG$ within the interior of the positive definite cone. Problem~\eqref{eq:LaplacianEstimation} can be solved using the alternating direction method of multipliers (ADMM). As a post-processing step, the minimizer can be fed into a minimum spanning tree algorithm to recover a tree graph. 

The third approach recovers $\bG$ by fitting data while ensuring the estimated $\bG$ is the reduced Laplacian of a tree graph via a mixed-integer linear program (MILP)~\cite{Taheri_MILP_2019}. To simplify notation, rewrite $\bG$ from \eqref{eq:G+B:G} as $\bG=\bA^\top\mcD_{\brho}\bA$ where $\brho$ collects the inverse line resistances. If $\bA$ is known, vector $\brho$ can be found via an LS fit. To this end, vectorize both sides in the RHS of \eqref{eq:m2} to get
\begin{equation}\label{eq:vectorize}
\bbp=\bH\brho+\bbn\quad\quad\text{where}\quad\quad\bH:=\tbV^\top\bA^\top*\bA^\top
\end{equation}
and $\bbp:=\vectorize(\tbP)$ and $\bbn:=-\vectorize(\bG\tbN)$. This follows from the property that $\vectorize (\bX\mcD_{\by} \bZ^\top)=(\bZ*\bX)\by$ where symbol $*$ denotes the Khatri-Rao matrix product. If $\bH$ is full column-rank, the LS estimate of $\brho$ is~\cite{Taheri_MILP_2019}
\[\brho^*_\bA=\arg\min_{\brho} \|\bbp-\bH\brho\|_2^2 = (\bH^\top\bH)^{-1}\bH^\top\bbp.\]
This problem alone solves the impedance estimation task. Different from topology identification or topology detection, impedance estimation amounts to the task of estimating line impedances assuming the topology to be known. References~\cite{sandia2,sandia3} provide impedance estimation solutions taking into account practical information, such as geographical/length constraints and line types a utility may be using. 

For topology identification, however, we also need to find $\bA$. Let us plug $\brho^*_\bA$ back into the LS fitting cost to get
\begin{equation}\label{eq:MILP1}
\|\bbp-\bH\brho^*_\bA\|_2^2=\|\bbp\|_2^2-\bbp^\top\bH(\bH^\top\bH)^{-1}\bH^\top\bbp.
\end{equation}
and minimize it over $\bH$ that depends on $\bA$. Recall each row of $\bA$ corresponds to an energized line. Suppose the operator has a library of $L>N$ lines that could be possibly energized. Candidate lines are encoded as rows of an $L\times N$ reduced branch-bus incidence matrix $\bbA$. Introduce binary variable $\bw\in\{0,1\}^L$ whose $\ell$-th entry indicates whether candidate line $\ell$ is energized. Finding $\bA$ now amounts to finding $\bw$. Minimizing \eqref{eq:MILP1} over $\bA$ can be shown to be equivalent to~\cite{Taheri_MILP_2019}
\begin{align}\label{eq:MILP2}
\min_{\bw,\bx}~&~\bbp^\top\mcD_\bw\bbp+\bbp^\top\mcD_\bw\bx\\
\textrm{s.to}~&~(\bC-\mcD_\bw)\bx=\mcD_\bw\bbp\nonumber\\
~&~\text{radiality constraints on }\bw\in\{0,1\}^L.\nonumber
\end{align}
where matrix $\bC$ depends on the line library and voltage data as $\bC:=[(\bbA\tbV\tbV^\top\bbA^\top)\circ(\bbA\bbA^\top) + \bI_{L}]^{-1}$ and $\circ$ denotes the Hadamard (entrywise) matrix product. Problem~\eqref{eq:MILP2} can be reformulated as a MILP upon applying McCormick linearization on the products between continuous and binary variables $\mcD_\bw\bx$. The formulation is correct as long as $\bb$ yields a connected radial graph, which can be enforced via the radiality constraints discussed in~\cite[Cor.~1]{STKSL22}. Although presented for smart meter data under probing, the approach is applicable to synchronized data setups \emph{S1)} and \emph{S2)}.

To summarize setup \emph{M2)}, if data exhibit high signal-to-noise ratios (SNR), the graph approach is the method of choice as it ensures exact topology recovery at polynomial complexity. For moderate to low SNR, the operator may have to select between the convex heuristic of \eqref{eq:LaplacianEstimation} and the MILP depending on the complexity/accuracy trade-off they can afford.

\subsection{Setup M3) Partial Data}\label{subsec:m3}
A more realistic setup is that non-synchronized data (voltages and injections alike) are collected only at a subset of buses $\mcS$. Given data $(\bV_\mcS,\bP_\mcS,\bQ_\mcS)$, the correlation-based approach of \eqref{eq:corr}--\eqref{eq:corrmn} can be used to compute the topology submatrices $(\bR_{\mcS\mcS},\bX_{\mcS\mcS})$. It is worth noting that correlation-based estimation of $(\bR_{\mcS\mcS},\bX_{\mcS\mcS})$ is viable even if injections in $\mcS$ are correlated (non-diagonal covariance matrices), as long as they remain uncorrelated from unobserved injections in $\bmcS=\mcN\setminus \mcS$ \cite{park2018learning}. This is relevant when $\mcS$ contains buses with solar generations that can exhibit strong correlation, whereas $\bmcS$ contains load buses alone that can be assumed uncorrelated over shorter intervals. References~\cite{park2018exact,Park_TCONES_2020,Cavraro_LCSS_2018} provide provable algorithms to derive the correct radial topology as long as all unobserved nodes in set $\bmcS$ have degree greater than $2$. This implies that all leaf nodes are in $\mcS$. For all leaf node pairs $\{m,n\}$, a well-known quantity termed \emph{effective resistance} can hence be computed as~\cite{Dorfler_2013_TCS} 
\[\bR_\text{eff}(m,n) = \bR(m,m)+\bR(n,n)-2\bR(m,n).\]

\begin{figure}[ht]
\centering	
\includegraphics[width=0.488\textwidth]{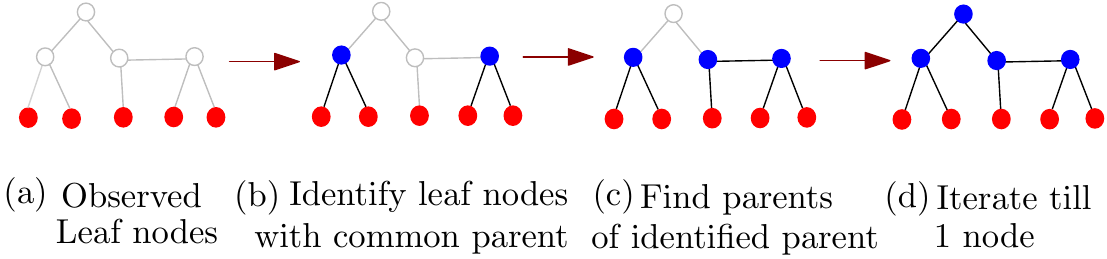}
\caption{Red nodes represent metered leaf nodes; rest of the nodes are non-metered. This figure illustrates the steps of iterative topology identification under setup \emph{M3)}, using effective resistances as a graph metric. The discovered nodes are colored blue, in each stage}
\label{fig:recursive}
\end{figure}

\begin{figure}[hbt]
\centering	
\includegraphics[width=0.5\columnwidth]{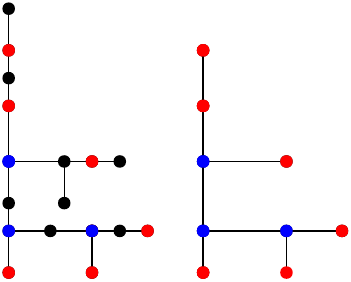}
\caption{Actual (left) versus recovered (right) grid topology by metering buses shown in red. The recovered topology remains radial. It covers all metered buses, as well as non-metered buses of degree at least three shown in blue. Buses in black cannot be recovered. All pairwise resistances on the recovered topology agree in value with the pairwise resistances on the original topology; this is important for voltage control applications.}
\label{fig:graph2}
\end{figure}

For radial grids and using \eqref{R_add}, effective resistance can be proven to be additive along graph paths, meaning $\bR_{\text{eff}}(m,k) =\bR_{\text{eff}}(m,n) +\bR_{\text{eff}}(n,k)$ for any node $k\in \mathcal{P}_{m,n}$. Additive nodal distances enable the construction of the radial graph using an iterative procedure termed \emph{recursive grouping (RG)}; see~\cite{choi2011learning,Pengwah2022SJ}. Specifically in our grid setting, at the first iteration of RG, each leaf node pair $(m,n)$ with a common parent can be identified by checking if $\bR_{\text{eff}}(m,k)-\bR_{\text{eff}}(n,k)$ is fixed for all other nodes $k$. If so, the parent of each pair $(m,n)$ is discovered, and the effective resistances between the parents and other observed nodes are derived. During the next iteration, the discovered parents are treated as leaf nodes, and the algorithm uses their effective resistances to identify common parents and so on. A schematic is provided in Fig.~\ref{fig:recursive}. While leaf node observations are sufficient, the algorithm is able to consider the case where intermediate nodes are also observed. Furthermore, if some non-metered node $k\in\bmcS$ has degree $1$ (leaf) or $2$, the algorithm is able to estimate a Kron-reduced graph without node $k$; see~\cite{Cavraro_LCSS_2018}. The recovered graph remains radial, maintains the relative ordering of all other nodes, and their pairwise effective resistances are found correctly; see Fig.~\ref{fig:graph2}. Reference~\cite{7981061-R2} proposes a graph processing algorithm to identify the grid topology from submatrices $\bR_{\mcS\mcS}$ or $\bX_{\mcS\mcS}$ building on \emph{Pr\"{u}fer sequences}, which resemble recursive grouping.

\subsection{Setup M4) Voltage-only Data}\label{subsec:m4}
A fourth setup that is of practical interest is when the operator accesses voltage data alone. This could the case because the measurement is taken at the photovoltaic panel, and not at the meter of a home. In this case, the meter may be measuring the solar generation, but not the net injection. Another scenario is that the measurement is actually a smart meter recording the total (re)active injection for a house. From the voltage and the injection data, an operator can extrapolate to estimate the voltage at the pole transformer if the secondary line characteristics are known. In this way, the operator can estimate the voltage at the distribution transformer even if such transformer serves 5-20 houses and not all of them are instrumented with smart meters. Given data $\bV$, the operator can only compute the centered correlation or covariance matrix
\begin{equation}\label{eq:corr2}
\bSigma_{vv}=\bR\bSigma_{pp}\bR+\bX\bSigma_{pq}^\top\bR + \bR\bSigma_{pq}\bX + \bX\bSigma_{qq}\bX+\bSigma_{nn}.
\end{equation}
We next review three approaches that rely on $\bSigma_{vv}$ to identify topologies using graph algorithms.

\emph{Signatures in Inverse Voltage Covariance.} Reference~\cite{Bolognani_CDC_2013} uses $\bSigma_{vv}$ to unveil a radial grid topology under the assumptions that: \emph{a1)} lines have identical resistance-to-reactance ratios ($r_\ell/x_\ell=\gamma$ for all $\ell$) so that $\bX=\gamma\bR$ from \eqref{eq:G+B}; \emph{a2)} power injections are independent across buses (diagonal $\bSigma_{pp}$, $\bSigma_{qq}$, and $\bSigma_{pq}$); and \emph{a3)} the measurement/modeling noise term $\bSigma_{nn}$ is negligible compared to other terms in \eqref{eq:corr2}. Under these assumptions, the voltage covariance turns out to be $\bSigma_{vv}=\bX\bSigma\bX$ for a diagonal matrix $\bSigma$ with positive diagonal entries that depend on injection covariances and $\gamma$. Due to this, the inverse covariance (a.k.a. \emph{precision matrix} in statistics) becomes $\bSigma_{vv}^{-1} = \bB\bSigma^{-1}\bB$. Thus, the precision matrix $\bSigma_{vv}^{-1}$ features the weighted structure of a \emph{squared} graph Laplacian of a radial network. It can be shown that the $(m,n)$-th entry of $\bSigma_{vv}^{-1}$ is positive if $m=n$ or buses $(m,n)$ have a common neighbor; negative if buses $(m,n)$ are neighbors; and zero, otherwise. In other words, negative valued entries in $\bSigma_{vv}^{-1}$ mark true edges of the radial grid. Reference~\cite{Bolognani_CDC_2013} suggests inverting the sample estimate of $\bSigma_{vv}$ and feeding its negative off-diagonal entries into a minimum-spanning algorithm to recover the grid topology. Numerical tests on a grid violating the common resistance-to-reactance ratio assumption demonstrate satisfactory topology identifiability probabilities for large $T$. 

Unfortunately, the precision matrix $\bSigma^{-1}_{vv}$ can be estimated as the inverse of the sample covariance $\hbSigma_{vv}=\frac{1}{T}\sum_{t=1}^T\bv_t\bv_t^\top$ only if $T\geq N$ and the grid undergoes sufficiently rich conditions so that $\hbSigma_{vv}$ is invertible. Nevertheless, if $\bSigma^{-1}_{vv}$ is known to be sparse, it can be recovered even from singular $\hbSigma_{vv}$ with proper regularization. To be more specific, if voltage data are modeled as Gaussian random vectors that are independent across time,\footnote{Voltages are approximately linear combinations of powers per \eqref{eq:ldf0:v}. If vectors $(\bp,\bq)$ are modeled as random, then $\bv$ tends to be Gaussian upon invoking a particular rendition of the central limit theorem~\cite[Remark~1]{Cavraro_TSG_2019}.} the precision matrix can be estimated via a sparsity-promoting maximum likelihood approach as the minimizer of the convex problem~\cite{tibshirani1996regression}
\begin{equation}\label{graph_lasso}
\bSigma^{-1}_{vv} = \arg \min_{\bS}\trace( \bS\hbSigma_{vv}) +\lambda\|\bS\|_1-\log|\bS|.
\end{equation}
This so-termed \emph{Graphical Lasso} method can estimate sparse precision matrices even with small $T$; see \cite{wainwright2008graphical}. Its computational complexity scales as $O(N^3)$~\cite{tibshirani1996regression}. Instead of using Graphical Lasso, each row of the precision matrix can also be estimated by $l_1$-norm regularized regression through the \emph{neighborhood Lasso}~\cite{meinshausen2006high}, where the voltage sequence at each bus is projected on all other buses as described in~\cite{Liao_PESGM_2016}. 

Reference~\cite{Deka_PSCC_2016} generalized the previous approach to show that when working with voltage phasor data, the assumption of constant $r/x$ ratio at all lines can be waived. In this setting, one can compute covariances between voltage magnitudes and angles and invert the related matrix to get
\begin{equation*}
\bSigma^{-1}_{(v,\theta)} = \begin{bmatrix} \bJ_{1} &\bJ_{2}\\\bJ^\top_{2} &\bJ_{3}\end{bmatrix}~~~\text{where}~~~\bSigma_{(v,\theta)} := \begin{bmatrix} \bSigma_{vv} &\bSigma_{v\theta}\\\bSigma_{\theta v} &\bSigma_{\theta\theta}\end{bmatrix}.
\end{equation*}
It can be shown that edge $(m,n)$ exists if and only if $\bJ_1(m,n) + \bJ_3(m,n) <0$, and thus, edges can be identified by negative entries in the inverse covariance. Theoretically, non-adjacent nodes $(m,n)$ with common neighbors (two-hop neighbors) have a strictly positive non-zero entry at $\bJ_1(m,n) + \bJ_3(m,n)$. Entries corresponding to node pairs being more than two hops away are $0$. Related works \cite{Liao_PESGM_2016,liao2018urban} showed that for realistic distribution grids, the values of $\bJ_1(m,n) + \bJ_3(m,n)$ at two-hop nodes are much smaller in magnitude than entries corresponding to adjacent nodes. In other words, thresholding the magnitude of entries in the inverse covariance enables topology detection as
 \begin{equation}\label{simple_CI}
(m,n) \in \mcL \text{~if and only if~} |\bJ_1(m,n)| + |\bJ_3(m,n)|>\lambda 
\end{equation} 
where $\lambda$ is a user-specified threshold \cite{Liao_PESGM_2016,liao2018urban}. The inverse covariance matrix $\bSigma^{-1}_{(v,\theta)}$ bears extra properties that can be exploited to recover grid topology without resorting to sign-based rules as described next.

\begin{figure}[!ht]
\centering
\includegraphics[width=0.11\textwidth]{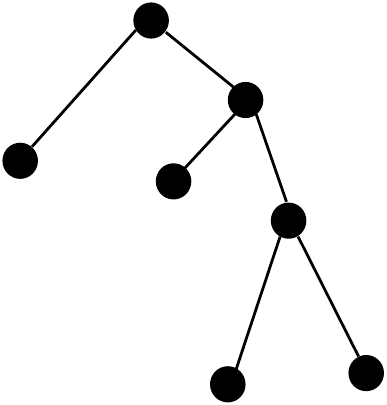}~~~~~~~~
\includegraphics[width=0.11\textwidth]{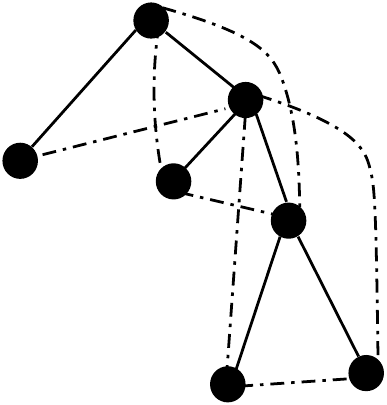}
\caption{\emph{Left:} Nodes in distribution grid $\mcG_o$. \emph{Right:} Graphical model or conditional independence (CI) graph $\mcG_{\text{CI}}$ for voltages. It includes tree edges (solid lines) and edges between 2-hop neighbors (dotted lines). Dotted lines need to be removed from CI graph to obtain grid graph.}
\label{fig:CI_voltage}
\end{figure}

\emph{Conditional Independence Test for Voltages.} If voltage data is modeled or approximated as Gaussian random variables, then zero entries in $\bSigma^{-1}_{(v,\theta)}$ imply that the corresponding variables are conditionally independent; see \cite{wainwright2008graphical,Deka_PSCC_2016} and references therein. Mathematically, $\bSigma^{-1}_{(v,\theta)}(m,n) = 0$ implies that $\mathbb{P}(v_m,v_n|u_k~\forall k\notin \{m,n\}) = \mathbb{P}(v_m|u_k~\forall k\notin \{m,n\})\mathbb{P}(v_m|u_k~\forall k\notin \{m,n\})$. Conditional independence tests can be accomplished using the inverse covariance matrix or computing mutual information.
The graph with edges for every pair of nodes that are not conditionally independent is called the \emph{conditional independence graph} or \emph{graphical model}~\cite{wainwright2008graphical}. It follows from the discussion in the previous paragraph that both adjacent nodes (one-hop neighbors), as well as two-hop neighbors in the grid $\mcG_o$, have edges in the conditional independence graph $\mcG_{\text{CI}}$, as shown in Fig.~\ref{fig:CI_voltage}.

Once the conditional independence graph of voltages has been derived, a conditional independence test \cite{Deka_PSCC_2016} can be used to distinguish true from spurious grid edges in $\mcG_{\text{CI}}$. Specifically, if voltages at two nodes $(m,n)$ are able to render voltages at some other nodes $(k,l)$ conditionally independent, i.e.,
$p(v_k,v_l|u_m,u_n) = p(v_k|u_m,u_n)p(v_l|u_m,u_n)$, then $(m,n) \in \mcL$. Two-hop neighbors do not satisfy this rule, and hence, can be identified and ruled out, while edges between non-leaf nodes (internal nodes) are identified. Once that is completed, edges between leaf nodes and their parents can be further determined. 

References~\cite{he2011dependency,Liao_MAPS_2015,Weng_TPS_2017} present a similar but simpler conditional independence test for topology recovery in realistic grids by showing that two-hop neighbors are almost conditionally independent; see \eqref{simple_CI}. As such, the true edges can be directly estimated by identifying node pairs $(m,n)$ for which their conditional dependence is above a threshold. A similar method using Markov random field (that characterizes conditional independence) is proposed and validated in \cite{Zhao_TSG_2020}. While performance is often subpar, the advantage of combinatorial learning algorithms is that they can be used to provide sample guarantees on performance~\cite{Deka_PSCC_2016,Weng_TPS_2017}. They can also be easily extended to multiphase setups or setups derived from grid dynamics; see Sections~\ref{subsec:dynamics}.

In recent works, inverse voltage covariances pertaining to partial data have been used for topology identification. These works extend the sparsity results on $\bSigma_{(v,\theta)}^{-1}$. The observer has access only to the submatrix of $\bSigma_{(v,\theta)}$ corresponding to set $\mcS$ and denoted as $\bSigma_{(v,\theta)}(\mcS,\mcS)$. If the non-metered buses comprising set $\bmcS$ are sufficiently less and non-adjacent, it has been shown in~\cite{doddi2020learning,anguluri2021grid} that \[[\bSigma_{(v,\theta)}(\mcS,\mcS)]^{-1}= \bSigma^{-1}_{(v,\theta)}(\mcS,\mcS) + \bL\]
where $\bSigma^{-1}_{(v,\theta)}(\mcS,\mcS)$ is a sub-matrix of the complete inverse covariance, and $\bL$ is a low-rank matrix. Upon computing $\bSigma_{(v,\theta)}(\mcS,\mcS)$ using partial voltage data, regularized convex optimization schemes (expanding on \eqref{graph_lasso}) can separate the sparse component $\bSigma^{-1}_{(v,\theta)}(\mcS,\mcS)$ from the low-rank component $\bL$. While the former provides the edges between observed nodes (set $\mcS$), the latter can be further decomposed to identify edges connected to $\bmcS$. 

\emph{Greedy Algorithms using Variance of Voltage Differences.} In addition to inverse voltage covariances, the statistic
\begin{equation}\label{eq:phi}
\phi_{mn}= \mathbb{E}\left[\left(v_m-v_n-\mathbb{E}[v_m-v_n)]\right)^2 \right]
\end{equation}
has been utilized for topology identification, oftentimes with better performance at low sample regimes~\cite{Deka_TCONES_2018,Deka_ECC_2016,Deka_TCONES_2020}. The statistic $\phi_{mn}$ is shown to be monotonically increasing along the paths on a radial grid. In other words, if path $(m,n)$ passes through node $k$, then both $\phi_{mk}$ and $\phi_{kn}$ are less than $\phi_{mn}$. Such increasing trends enable a spanning tree-based topology recovery algorithm: Give each node pair $(m,n)$ a weight of $\phi_{mn}$ and construct the minimum spanning tree from the set of weights greedily. This will coincide with the true radial topology. As no inverse matrix operation is needed, the performance of the spanning tree algorithm at small numbers of voltage samples is often improved over other voltage-only learning algorithms. We next discuss topology detection schemes and draw connections to the identification schemes presented thus far. 

\section{Topology Detection}\label{sec:detection}
\begin{figure}[hbt]
\centering
\includegraphics[width=0.5\columnwidth]{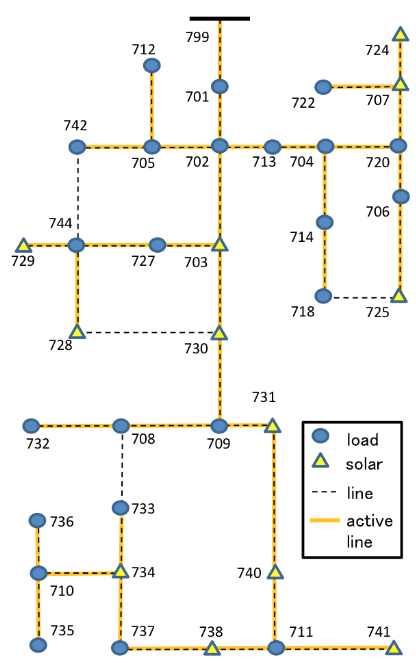}
\caption{An augmented IEEE 37-bus feeder benchmark from~\cite{Cavraro_TSG_2019}. Black dotted edges show line infrastructure. Lines shown in yellow are actually energized.}
\label{fig:detection}
\end{figure}

Compared to topology identification, topology detection is a relatively simpler task. The operator now already knows the line infrastructure, namely the impedance of each line as well as the buses incident to it. The operator does not know the status of all or a subset of lines. Given data, topology detection amounts to deciding which of the existing lines are energized (see Fig.~\ref{fig:detection}). This is useful during normal grid operation due to automated switching. Variants of topology detection schemes are also relevant after extreme weather events to inform grid restoration efforts for line outages and dispatch line crews. 

A general formulation parallels the modeling after \eqref{eq:MILP1} by including the now available information on lines. Often, this information is available as a library of known lines $\mcL_{\text{all}}$ with cardinality $L_{\text{all}}$, and known impedances. Out of these lines, the operator would like to select a set of lines $\mcL$ of size $N$ that are energized. Each candidate line can be associated with a row of an augmented $L_{\text{all}}\times N$ branch-bus incidence matrix $\bbA$. The status of line $\ell$ is indicated by a binary variable $w_\ell$. If vector $\bw$ and a diagonal matrix $\mcD_{\bw}$ collect all $L^1$ line statuses, the topology matrices of interest can be parameterized as
\begin{subequations}\label{eq:param}
\begin{align}
\bY(\bw) &= \bbA^\top \mcD_{\bw} \mcD_{\by} \bbA\\
\bR(\bw) &= (\bbA^\top \mcD_{\bw} \mcD_{\br}^{-1} \bbA)^{-1}\\
\bX(\bw) &= (\bbA^\top \mcD_{\bw} \mcD_{\bx}^{-1} \bbA)^{-1}.
\end{align}
\end{subequations}
Apparently, if the status of line $\ell$ is known, variable $w_\ell$ can be set to a fixed value. As with topology identification, topology detection methods are organized next based on the type of measurements used. Before proceeding with the methods, it is worth adding a note on topology detectability. Sections~\ref{sec:idsyncdata} and \ref{sec:idsmdata} discussed topology identifiability results, i.e., conditions under which the complete or reduced topology of a grid can be recovered from complete or partial data. Given that topology detection is easier than identification, detectability can be ensured with less strict conditions on the number and placement of meters as established in~\cite{Cavraro_LCSS_2020}.

\subsection{Topology Detection using Phasor Data}
Analogous to setup \emph{S1)} in Section \ref{subsec:S1}, if phasor data are collected at all buses, one can try detecting $\bw$ via an LS fit. Similar to \eqref{eq:vectorize}, upon vectorizing the data model in \eqref{eq:i=Yv} and exploiting the parameterization $\bY(\bw)$, we arrive at the mixed-integer quadratic problem (MIQP):
\begin{align}\label{eq:TopDetPMUProb}
\min_{\bw}~&~\|\vectorize(\bI) - \bH \bw\|_2^2\\
\textrm{s.to}~&~\text{radiality constraints on }\bw\in\{0,1\}^{L_{\text{all}}}.\nonumber
\end{align}
where $\bH:=\bU^\top\bbA^\top\mcD_{\by}*\bbA^\top$.

If the operator has to select among a few possible topologies $\{\bw_s\}_{s=1}^S$, one may evaluate which topology yields the smallest fitting error $\min_{s}\|\bI-\bY(\bw_s)\bU\|_F^2$, or the voltage reconstruction error $\min_{s}\|\bY^{-1}(\bw_s)\bI-\bU\|_F^2$ as suggested in~\cite{Arghandeh_PESGM_2015}.

\subsection{Topology Detection using Smart Meter Data}
The topology identification schemes presented earlier can be tailored for topology detection. For example, if the operator collects data using probing under setup \emph{M2)} [cf. Sec.~\ref{subsec:m2}], one can detect $\bw$ by modifying \eqref{eq:LaplacianEstimation}; see~\cite{Cavraro_TCONES_2019}. The regularization terms may not be needed anymore, and the formulation becomes similar to the one in \eqref{eq:TopDetPMUProb}. To avoid mixed-integer programs, one may also pursue a box relaxation on $\bw$.

Adjusting identification algorithms towards topology detection is straightforward for schemes under partial data setups \emph{M2)}--\emph{M3)} described in Sections \ref{subsec:m2} and \ref{subsec:m3}. As the set $\mcL_{\text{all}}$ of possible lines is known, edge discovery in the combinatorial graph algorithms is restricted to $\mcL_{\text{all}}$ instead of considering all node pairs. For example, in the recursive grouping-based scheme in \emph{M2)} and \emph{M3)}, known edges can be included beforehand and only permissible parent-child pairs are discovered. 

Similarly, \emph{detection algorithms for voltage-only data} can be designed by modifying the signature or conditional independence (CI)-based tests on inverse voltage covariances, discussed under setup \emph{M4)} in Sec.~\ref{subsec:m4}, to be restricted to $\mcL_{\text{all}}$. In an alternate approach, the dependence of permissible edges $\bw$ on $\bSigma^{-1}_{vv}$ can be directly parameterized through \eqref{eq:param} and \eqref{eq:corr2}. Similar to \eqref{graph_lasso}, given the sample voltage covariance $\hbSigma_{vv}$, a maximum likelihood detector is put forth in~\cite{Cavraro_TSG_2019} to select $N$ edges from $\mcL_{\text{all}}$:
\begin{align*}
\min_{\bw}&~\trace\left(\bSigma_{vv}^{-1}(\bw)\hbSigma_{vv}\right) -\log |\bSigma^{-1}_{vv}(\bw)|\\
\textrm{s.to}&~\bw\in [0,1]^{L_{\text{all}}},~~\bone^\top \bw=N.
\end{align*}
Despite the box relaxation on $\bw$, the objective function is non-convex, and a stationary point is found through a projected gradient descent algorithm. In a related work, \cite{Sharon_ISGT_2012} uses a maximum likelihood framework for Gaussian voltage distributions to identify topology changes from a set of permissible edges.

\begin{figure}[hbt]
\centering
\hspace*{\fill}
\includegraphics[width=0.15\textwidth]{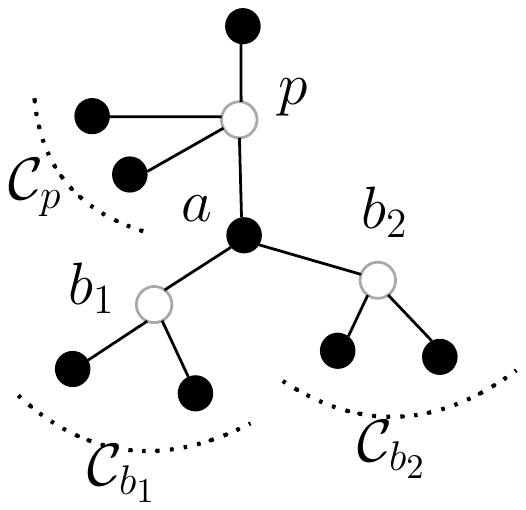}\hfill
\includegraphics[width=0.15\textwidth]{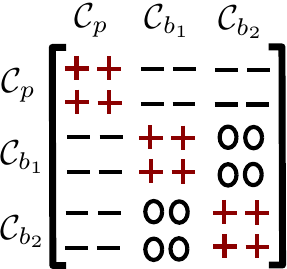}
\hspace*{\fill}
\caption{\emph{Left:} A radial tree with observed nodes colored solid black. $C_p$, $C_{b_1}$ and $C_{b_2}$ are children of nodes $p$, $b_1$ and $b_2$ respectively. $\phi_{mn}$, defined in \eqref{eq:phi}, increases as node $n$ is moved away from $m$, along any path, enabling topology learning if all nodal voltages are observed. \emph{Right:} The sign of entries for matrix $[\phi_{k_1a} +\phi_{k_2a}-\phi_{k_1k_2}]$ for $k_1,k_2 \in C_{b_1}\cup C_{b_2}\cup C_p$ is based on the relative location of $k_1, k_2$ with respect to $a$. This enables separating observed nodes into groups for topology detection.}\label{fig:missing}
\end{figure}

\textit{Detection with Partial Voltage and Injection Statistics.} The variances of voltage differences $\phi_{mn}$, defined in \eqref{eq:phi} under \emph{M4)}, can be used to detect topologies if available for all node pairs. Additionally, it was shown in \cite{Deka_ECC_2016,Deka_TCONES_2020,Deka_SmartGridComm_2016} that $\bphi$'s carry more information:
\renewcommand{\theenumi}{\emph{p\arabic{enumi}}}
\begin{enumerate}
\item The statistic $\phi$ has additional ordering for triples of nodes $(m,n,k)$, that is, the sign of $\phi_{mn}-\phi_{mk}-\phi_{nk}$ depends on the relative position of $(m,n,k)$ on the grid; see Fig.~\ref{fig:missing}.
\item The value of $\phi_{mn}$ at adjacent nodes $(m,n)$ with parent node $m$, depends on the admittance of edge $(m,n)$ and injections at all descendant nodes of $n$.\end{enumerate}

These ordering properties enable topology detection, where voltage and injection statistics can be computed over a subset of buses $\mcS$. Under partial data, reference~\cite{Deka_TCONES_2020} first creates a spanning tree for observed nodes using available $\phi$'s, and then inserts missing nodes iteratively by checking for properties \emph{p1)}--\emph{p2)}. The algorithm is able to detect the true topology provided each unobserved node has degree greater than $2$. This is similar in theme to the learning algorithms proposed under recursive grouping in \emph{M3)}, with the crucial distinction that the algorithm has access to injection statistics and line impedances in $\mcL_{\text{all}}$, but not injection samples. The latter means that $\bSigma_{vp}$ and hence $\bR$ cannot be estimated. 

If we are confined to a few candidate topology vectors $\{\bw_s\}_{s=1}^S$, reference~\cite{He2021hybrid} suggests a model validation process, wherein the power flow solution obtained using each candidate topology is compared with the collected grid data to determine the operating topology. A similar workflow is pursued in \cite{9852281} using Bayesian networks.

\subsection{Topology Detection using Line Flow Data}
Previously, we assumed that sensors have been installed only at buses. However, some works cope with topology detection using also line flow measurements. The work~\cite{Farajollahi_TSG_2020} considers a setup where only line current measurements and load (pseudo) measurements are available. Thanks to a power flow linearization, they cast topology detection as a MILP to minimize the mismatch between measured and computed line currents. Reference~\cite{sevlian2015distribution} exploits the fact that in radial networks and upon neglecting power losses, the power flow on a line $(m,n)$ equals the sum of powers of all downstream buses. Given line flow and voltage data at a few locations along with (pseudo) measurements of loads, reference~\cite{Tian_TPS_2016} casts topology detection as a mixed-integer quadratic program.

\section{Extensions}\label{sec:extensions}
The schemes discussed so far refer to balanced or single-phase and radial grids, presumed steady-state operation, and relied on readings of electric and particularly nodal quantities (voltages, injections). We next briefly discuss topology learning schemes that waive one or more of these assumptions. 

\subsection{Learning Loopy/Meshed Topology}\label{subsec:loopy}
While distribution grids are typically operated under radial topologies, loopy or meshed topologies can be found at heavily loaded urban distribution grids~\cite{NY}. Despite complications related to protection schemes, meshed topologies could advance reliability, especially in micro-grid settings. In this context, the inverse voltage covariance-based topology identification methods discussed in Section \ref{subsec:m4}, can also learn meshed grids under realistic restrictions on the minimum size of loops/cycles, as shown in \cite{liao2018urban,he2011dependency,Deka_TSG_2020}. Further, these methods also extend to meshed grids with ZIBs (zero-injection buses), where the voltage covariance matrix is non-invertible. In that case, reference~\cite{Deka_TCONES_2021} first identifies ZIBs and their neighbors using flow conservation-based voltage regression, and then, discovers the remaining topology using properties of the inverse voltage covariance at non-ZIBs.

\subsection{Multiphase Topology Learning}\label{subsec:multiphase}
The linearized power flow model can be extended to multiphase settings under similar assumptions of small angle differences (per phase) for neighboring grid buses~\cite{VKZG16,gan2014convex}. The primary complexity under a multiphase model arises from the fact that voltages and injections at each bus are not complex scalar quantities anymore, and each line admittance is a complex-valued $3\times 3$ or $2\times 2$ matrix with off-diagonal entries. Furthermore, the phase assignment of each bus may be unknown. While an extensive discussion on \emph{phase identification} falls outside the scope of this article, we note that methods based on voltage clustering and local voltage-based distances have been employed to assign phases~\cite{blakely2019spectral,wang2016phase,olivier2018phase, bariya2021guaranteed}. If phases are correctly identified or known, the variance of nodal voltage differences (see the single-phase case under setup \emph{M4)} in Sec.~\ref{subsec:m4}) can be used as edge weights to identify the topology as a spanning tree. Crucially, even without phase identification, the conditional independence structure of three-phase nodal voltages factorizes exactly as described for the single-phase case in Sec.~\ref{subsec:m4} and Fig.~\ref{fig:CI_voltage}. Conditional independence tests for three-phase voltages have been proposed for topology identification in~\cite{Deka_TPS_2020,liao2019unbalanced}, using properties akin to the single-phase setting. 

A mixed-integer linear program (MILP) formulation for three-phase topology detection is presented in~\cite{gandluru2019joint}, that can be categorized as a three-phase extension of the partial observability setup \emph{M2)} in Sec.~\ref{subsec:m2}. The MILP objective here minimizes the weighted error in three-phase linearized power flow equations involving state variables, smart meter data, and pseudo-measurements. Reference~\cite{Li_JEPES_2021} considers a regime where three-phase voltages and injections at leaf nodes alone are available, similar to the single-phase setup \emph{M2)} in Sec.~\ref{subsec:m3}. The topology identification uses recursive grouping but is applied to efficiently estimate in-phase effective impedances. 

\subsection{Learning with Multiple Feeders} 
Distribution networks are typically operated as a union of radial sub-networks, each fed by a different substation or substation transformer referred to as source point. However, only a few works can cope with topology learning if the knowledge of which bus is fed by what source point is not known in advance~\cite{Deka_TCONES_2018,Cavraro_TSG_2019}. Such information might be lacking because utilities may have limited information about secondary circuits, and loads are oftentimes transferred between feeders to alleviate overloads, or during restoration after extreme weather events. Hence, before performing any topology processing algorithm, one should cluster in advance buses based on the source point to which they are connected. The main idea here is that the power flowing through the source point matches, up to the power loss, to the sum of the served loads~\cite{Cavraro_TSG_2020,arya}, and that voltages across feeders can be approximated as uncorrelated~\cite{Deka_TCONES_2018,Cavraro_TSG_2020}.

\subsection{Learning with Grid Dynamics}\label{subsec:dynamics}
All previous works considered the steady-state operation of the distribution grid. In a fast time-scale setting, reference~\cite{bariya2018data} compares the voltage fluctuations experienced at each bus after events, such as transformer tap or capacitor bank switching. It is observed that the closer two buses are, the more similar their responses to such events are; hence this can be used as a metric for topology identification. In fact, as grid reconfiguration itself induces changes in voltage profiles, researchers have devised methods to detect topology changes by comparing voltage data with a library of signatures \cite{Cavraro_ISGT_2015,Cavraro_PESGM_2015,Cavraro_TPS_2018}, or using a subspace perturbation model~\cite{zhou2019power}.

\subsection{Topology Processing with Non-Electric Data}
So far, we have discussed how grid topologies can be inferred by processing electric grid data. Nonetheless, grid topologies can also be estimated using non-electric data. For example, reference~\cite{Erseghe_TSP_2013} uses power line communication (PLC) signals to identify neighboring nodes on a radial grid. It first measures the communication delay $t_{nm}$ for a PLC signal to be propagated between any pair $(n,m)$ of buses. Based on pairwise delays, two buses $(n,m)$ are declared as neighbors if $t_{nm}<t_{nk}+t_{km}$ for any other bus $k\notin\{n,m\}$. The key principle here is that the communication delay scales with the line length.

\section{Topics Not Covered and Future Directions}\label{sec:topics}

This tutorial has not discussed the considerable amount of work on topology learning in \emph{transmission systems}. For example, the classical task of generalized state estimation (GSE) aims at jointly estimating the system state along with its topology~\cite{AburExpositoBook,Kekatos_NAPS_2012}. Beyond GSE, there have been efforts towards learning the topology of transmission grids using various types of data. For example, reference~\cite{Tirja2} uses voltage angle data to recover the network topology using an optimization scheme similar to \eqref{graph_lasso}. Reference~\cite{Kekatos_TSG_2016} exploits the structure of wholesale electricity markets to infer the topology of power transmission systems using locational marginal prices. References~\cite{Talukdar_ACC_2017,Talukdar_Automatica_2020} attempt to recover network topologies from voltage angle data exploiting the structure of swing dynamics. While phase angle samples now are temporally correlated, conditional independence-based tests, similar to ones in Section~\ref{subsec:m4}, can identify neighboring buses even in meshed grids. The key difference here is that instead of inverse covariance, the sparsity and signs of the inverse \emph{power spectral density} matrix are exploited; power spectral density is the Fourier transform of auto-correlation $\mathbb{E}[\btheta_t\btheta^T_{t+h}]$. Conditional independence-based tests on inverse power spectral density have been extended for consistent identification in grids with partial observability~\cite{Talukdar_ACC_2018}, as well.

Focusing again on distribution grids, there is extensive literature on other topology-processing topics that were only briefly discussed or not touched upon at all here. For example, the concept of GSE has been extended to distribution grids. As with topology detection, line impedances are assumed to be known. It is performed by solving weighted least-squares or least-absolute-value (WLAV) problems over continuous (system states) and binary (line or breaker statuses) variables~\cite{9748966}. If the measurement set does not ensure identifiability, pseudo-measurements are oftentimes introduced~\cite{Korres2012GSE}. Reference~\cite{Soltani2022GSE} casts GSE as a mixed-integer quadratic program for unbalanced multiphase feeders with possibly missing/bad data. 

In addition, there is a growing literature on model-agnostic and supervised topology processing techniques, that has not been covered herein. Such techniques typically proceed in two stages. During the training stage, they collect labeled grid data~\cite{Dua21}. Labels may correspond to known outages, particular topologies, or fault locations. Such labeled data can be used either as a library of recorded events; or alternatively, used to train a machine learning model (such as a deep neural network, DNN) in a supervised manner for classification or regression; see e.g., \cite{Azimian22} and references therein. During operation, the operator uses the current data and the learned model to infer grid topology information. Most of the methods reviewed in this article presumed customers are Wye-connected and hence may need modifications to deal with Delta or line-to-line connections in general. Fault localization has been an area of active research in power transmission systems, with increasing relevance in distribution grids due to several grid-initiated wildfires in recent years. Detecting and identifying such faults with limited phasor data is a crucial research direction. 

Moving forward, there are several research pathways in learning distribution grid topologies that warrant further investigation. Research directions for distribution grid topology learning that have not been studied as extensively include but are not limited to nonlinear power flow models or setups where grid data exhibit high correlation (e.g., due to solar injections) or suffer from low signal-to-noise ratios. Moreover, learning microgrid topologies is a pertinent task given their potentially mobile and ad-hoc nature. Unveiling meshed topologies may become increasingly more relevant for densely populated urban distribution grids. On the methodological front, tensor processing techniques may prove fruitful towards unveiling grid topologies.

New metering technologies and the aggregation of large amounts of heterogeneous data offer unprecedented opportunities for topology learning. Collecting synchrophasor data at higher-order harmonics in addition to the phasors at the fundamental frequency may offer neat opportunities for grid topology learning. New meters that can record phasors accurately at sub-cycle levels (termed point-on-wave measurements) have led to growing interest in using them to identify topologies, estimate impedances, and deal with affiliated tasks, despite limited placement. Localizing hidden sources of oscillations has also become pertinent given new types of loads and inverter-interfaced distributed energy resources (DERs). PMUs or sensors with synchrophasor functionality installed in DERs provide data on the grid edge. Smart inverters and other smart grid devices (e.g., smart transformers, remotely-controlled voltage regulators and capacitor banks, D-STATCOMs) offer additional opportunities to collect data and estimate voltage sensitivities upon actuation and probing. Flow measurements from transformers and line sensors add to the abundance of data. Non-electric data such as GPS (global positioning system) data from AMI as well as data from drones surveying lines bring in a spatial dimension to the topology learning task at hand. Thermal cameras sensing the temperature of lines from afar can potentially infer whether a line is energized or not. Fusing data from heterogeneous sources, such as smart meters, phasor measurement units, line sensors, or geographical data from city maps and street layouts, aerial/drone images could be of practical interest.

As topology estimation methods are data-driven, targeted adversarial corruption of the measurement set can lead to errors in topology recovery. This is within the purview of data injection and jamming based attacks in power grids and warrants future research. Meter placement to improve topology estimation \cite{Dua21} and robustness is a potential area of future work. Topology identification and detection tools may be linked to feedback loops for data-driven grid control, where changes in structure or parameter values are seamlessly incorporated into evolving control logic. Studying the optimality and uncertainty quantification of such feedback presents an interesting research opportunity.

\section{Conclusion}\label{sec:conclusion}
This tutorial has presented a study of topology identification and detection methods in radial distribution grids. Focusing on analytical methods, we have compared and contrasted topology learning algorithms based on measurement type, the extent of observability, and the solution approach proposed. Under realistic low-observability regimes, we have presented explanatory discussions on the use of power flow physics and topological constraints in different algorithms to ensure topology recovery as well as discussed limitations on exact estimation. Connections to learning schemes for meshed, unbalanced grids, and related problems in change detection and physics-agnostic learning methods have been elucidated. Finally, we have commented on potential directions where the current state-of-the-art methods may be extended in the future. 

\balance

\bibliographystyle{IEEEtranS}
\bibliography{arxiv_v2}

\begin{IEEEbiography}[{\includegraphics[width=1in,height=1.25in,clip,keepaspectratio]{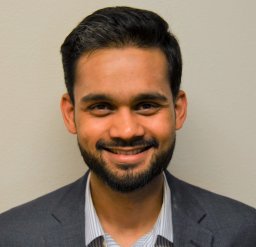}}]{Deepjyoti Deka} (SM'20) is a staff scientist in the Theoretical Division at Los Alamos National Laboratory (LANL), where he was previously a postdoctoral research associate at the Center for Nonlinear Studies. At LANL, Dr. Deka serves as a PI/co-PI for DOE and LDRD projects on machine learning in power systems, interdependent networks, and in cyber-physical security. Before joining the laboratory, he received his M.S. and Ph.D. degrees in electrical and computer engineering (ECE) from the University of Texas, Austin, TX, USA, in 2011 and 2015, respectively. He completed his undergraduate degree in electronics and communication engineering (ECE) from IIT Guwahati, India, with an institute silver medal as the best outgoing student of the department in 2009. Dr. Deka is a senior member of IEEE and has served as an editor on IEEE Transactions on Smart Grid in 2020-2021.
\end{IEEEbiography}

\begin{IEEEbiography}[{\includegraphics[width=1in,height=1.25in,clip,keepaspectratio]{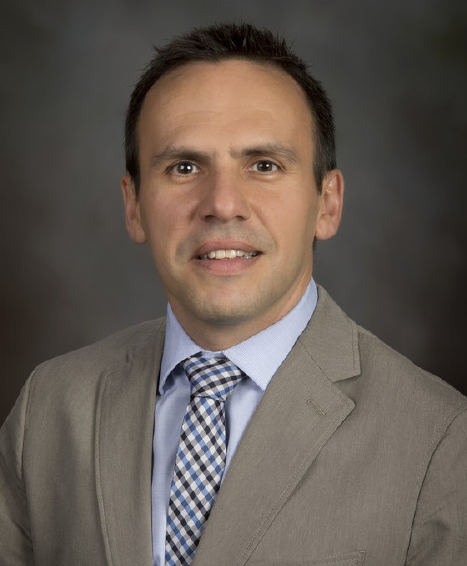}}] {Vassilis Kekatos} (SM'16) is an Associate Professor with the Bradley Dept. of ECE at Virginia Tech. He obtained his Diploma, M.Sc., and Ph.D. in computer science and engineering from the Univ. of Patras, Greece, in 2001, 2003, and 2007, respectively. He is a recipient of the NSF Career Award in 2018 and the Marie Curie Fellowship. He has been a research associate with the ECE Dept. at the Univ. of Minnesota, where he received the postdoctoral career development award (honorable mention). In 2014, he stayed with the Univ. of Texas at Austin and the Ohio State Univ. as a visiting researcher. His research focus is on optimization and learning for future energy systems. During 2015-2022, he served on the editorial board of the IEEE Trans. on Smart Grid.
\end{IEEEbiography}

\begin{IEEEbiography}[{\includegraphics[width=1in,height=1.25in,clip,keepaspectratio]{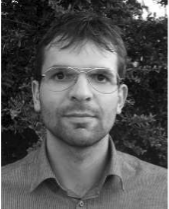}}] {Guido Cavraro} received the Ph.D. degree in Information Engineering from the University of Padova, Italy, in 2015. He was a visiting scholar at the California Institute for Energy and Environment (CIEE) at U.C. Berkeley in 2014. In 2015 and 2016, he was a postdoctoral associate at the Department of Information Engineering of the University of Padova. From 2016 to 2018, he was a postdoctoral associate with the Bradley Department of Electrical and Computer Engineering of Virginia Tech, USA. Currently, he is a Senior Researcher with the Power Systems Engineering Center at National Renewable Energy Laboratory, USA. His research interests include control, optimization, and estimation applied to power systems and smart grids.
\end{IEEEbiography}

\end{document}